\documentstyle[amssymb,amstex,12pt,a4wide]{article}
\newtheorem{theorem}{Theorem}
\newtheorem{lemma}{Lemma}
\newtheorem{corollary}{Corollary}

\newcommand{\tr}{\operatorname{Tr}}
\newcommand{\Po}{\operatorname{Po}}

\begin{document}
\title{Spectra of winner-take-all stochastic neural networks}
\author{Tomasz Schreiber\footnote{Mailing address: Tomasz Schreiber, Faculty of
 Mathematics \& Computer Science, Nicolaus Copernicus University,
 ul. Chopina 12 $\slash$ 18, 87-100 Toru\'n, Poland; tel.: (++48) (+56)
 6112951, fax: (++48) (+56) 6228979; e-mail: {\tt tomeks at mat.uni.torun.pl}},
 \footnote{
        Research supported by the Polish Minister of Science and Higher Education
        grant N N201 385234 (2008-2010)}\\
        Faculty of Mathematics \& Computer Science,\\
        Nicolaus Copernicus University, Toru\'n, Poland.}

\date{}
\maketitle

\setlength{\baselineskip}{1.2\baselineskip}

\paragraph{Abstract}
 {\it During the recent few years, in response to empirical findings suggesting scale-free self-organisation
       phenomena emerging in complex nervous systems at a mesoscale level, there has been significant
       search for suitable models and theoretical explanations in neuroscientific literature, see the recent
       survey by Bullmore \& Sporns (2009).  In Piekniewski \& Schreiber (2008) we have developed
       a simple and tractable mathematical model  shedding some light on a particular class of the
       afore-mentioned phenomena, namely on mesoscopic level self-organisation of functional brain
       networks under fMRI imaging, where we have achieved a high degree of agreement with
       existing empirical reports. Being addressed to the neuroscientific community, our work
       Piekniewski \& Schreiber (2008) relied on semi-rigorous study of information flow structure
       in a class of recurrent neural networks exhibiting asymptotic scale-free behaviour and
       admitting a description in terms of the so-called winner-take-all dynamics. The purpose of the
       present paper is to define and study these winner-take-all networks with full mathematical rigour
       in context of their asymptotic spectral properties, well known to be of interest for neuroscientific 
       community. Our main result is a limit theorem for spectra of the spike-flow graphs induced
       by the winner-take-all dynamics. We provide an explicit characterisation of the limit spectral
       measure expressed in terms of zeros of Bessel's J-function.} 

\paragraph{Keywords:} spectra of random scale-free graphs, winner-take-all dynamics,
                      neural networks.

\paragraph{MSC:} 60F05; 60K35; 15A52

\section{Introduction and motivations}
 Recent few years in the neuroscientific literature have been marked by a very successful interdisciplinary
 interaction between the study of large-scale phenomena in complex nervous systems and random graph
 theory, especially in context of the so-called scale-free networks considered a nearly classical subject by
 now, see e.g. Albert \& Barab\'asi (2002) or Chung \& Lu (2006) and Durett (2007) for a mathematical treatment. 
 Among a plethora of particular topics studied, the one in focus of our interest are the statistical properties
 of the so-called {\it functional brain networks} arising under fMRI imaging at mesoscale (usually understood
 as individual voxel level) where small world and scale-free self-organisation of activity correlations has
 been reported in empirical findings, see e.g. Bullmore \& Sporns (2009) for an extensive
 review and Egu\'iluz et al. (2005), Salvador et al. (2005), Cecchi et al. (2007) and van den Heuvel (2008)
 for presentation and discussion of experimental results. Certain heuristical non-rigorous considerations
 aimed at explaining these phenomena have been offered in Fraiman (2009) and Kitzblicher (2009) discussing
 very interesting analogies between crucial features of functional brain networks and Ising model
 at criticality. Up to our best knowledge, the first dedicated mathematical model shedding some light
 on the scale-free properties of mesoscopic brain functional networks is the simple spin glass type
 system introduced in Piekniewski \& Schreiber (2008) further extended and enhanced with a
 geometric ingredient in Piersa, Piekniewski \& Schreiber (2010) and standing in good agreement
 with empirical findings. The details and neuroscientific motivations of these models are far beyond
 the scope of the present mathematically oriented paper and we only provide a brief overview for
 completeness here, proceeding to well-defined rigorous problems as soon as possible. 

 The disordered system proposed in Piekniewski \& Schreiber (2008) models
 an asynchronous spiking neural network with the aim of analysing the structure
 of information flow in a class of recurrent neural nets. The model, bearing formal resemblance
 to the celebrated Sherrington-Kirkpatrick (1972) spin glass yet exhibiting quite
 different behaviour, consists of $N$ formal neurons $\varsigma_i,\; i=1,\ldots,N,$
 where the value $\varsigma_i \in \{0,1,2,\ldots\}$ represents the {\it charge} (activity level) 
 stored at $\varsigma_i.$ Initially each neuron stores some small fixed charge. The charge-conserving
 Kawasaki-style evolution of the system takes place by choosing at random subsequent pairs
 of numbers $i \neq j$ and trying to transfer a unit charge from $\varsigma_i$ to $\varsigma_j$
 \--- as soon as $\varsigma_i > 0$ such a
 trial is always successful if it decreases the energy of the system and is accepted
 with probability $\exp(-\beta \Delta H)$ and rejected with the complementary
 probability otherwise, where $\beta > 0$ is some positive {\it inverse temperature} parameter
 whereas the energy $H$ of the system is given by $H := \frac12 \sum_{i \neq j} w_{ij} 
 |\varsigma_i - \varsigma_j|$ with $w_{ij} = w_{ji}$ standing for i.i.d. standard Gaussian
 {\it connection weights}. In standard intuitive terms, the presence of a positive weight between 
 two neurons indicates that the system favours the agreement of their activity levels
 whereas a negative weight means that disagreement is preferred.  
 The object in the focus of our interest in Piekniewski \& Schreiber (2008)
 was the {\it spike flow-graph} or {\it charge-flow network} generated by this dynamics,  
 defined by ascribing to each edge $(ij)$ the multiplicity equal to the number 
 of charge transfers occurring along $(ij)$ in the course of (a long enough period of)
 the dynamics. This object gains a natural interpretation upon noting that edges with
 high multiplicities are those essential to the dynamics as designed to model
 the neural network's spiking activity, whereas the low multiplicity edges are
 only seldom used and could as well be removed from the network without 
 effectively affecting its evolution. In informal terms, the charge-flow graph
 represents the essential support of the system's effective dynamics, whence 
 our interest in this object. 

 In Piekniewski \& Schreiber (2008) we have performed a semi-rigorous analysis
 of the above model, based on extreme value theory methods, arguing that for
 $N$ large enough its ground state arises by putting the whole system charge
 into one {\it best} neuron (determined as a function of weights $w_{ij}$)
 and leaving all the remaining ones empty. Moreover,
 in low enough temperatures, the dynamics of such networks in large $N$ asymptotics
 is well approximated, in the sense made precise ibidem, by a much simpler
 {\it winner-take-all} (WTA) dynamics described in detail and rigour in
 Section \ref{MainRes} below. This observation allowed us to show in 
 Piekniewski \& Schreiber (2008) that asymptotically the charge-flow networks are scale-free
 with exponent $2,$ see ibidem as well as Piersa, Piekniewski \& Schreiber (2010), in agreement
 with the empirical findings as quoted above. We have also argued there that even though the
 spin glass model we propose may be regarded quite specific, its large scale behaviour and
 in particular its winner-take-all approximation is presumably universal for a large class of
 networks where each formal neuron represents a computational unit exhibiting some
 non-trivial internal structure and memory, for instance a group of biological or artificial neurons
 (see Piekniewski, 2007) whose internal state requires more complicated
 labeling than just $\{-1,+1\}$ as in the original Sherrington-Kirkpatrick
 model, whence the ${\Bbb N}$-valued labels in our model.

 The purpose of this paper is to complement the semi-rigorous developments of 
 Piek\-niew\-ski \& Schreiber (2008) by carrying out a fully rigorous mathematical study
 of the asymptotic structure of random charge-graphs generated by the winner-take-all
 dynamics described in full detail in Section \ref{MainRes} below. More precisely,
 we focus on spectral measures of these graphs as providing important information 
 about their underlying structure,  see e.g. Chapters 8 and 9 in Chung \& Lu (2006)
 for a  discussion of spectral aspects of scale-free graphs. 

\section{The model and main results}\label{MainRes}
To provide a formal description of the winner-take-all dynamics,
consider the set $\{ 1,\ldots,n \}$ of network vertices, each vertex identified
with its rank between $1$ and $n.$ Initially are $m = \lfloor \alpha n \rfloor,\;
\alpha \in {\Bbb R}_+,$ {\it units of charge} present in the system,
with each unit stored in a vertex chosen uniformly by random,
independently of other units. The system evolves thereupon
according to the following sequential {\it winner-take-all} (WTA)
dynamics, with $\sigma_i$ standing for the current charge stored at $i.$
\begin{description}
 \item {\bf (WTA)} Choose uniformly by random a source vertex $i \in \{1,\ldots,n\}$
                   and, independently, a target vertex $j \in \{1,\ldots,n\}.$
                   \begin{itemize}
                    \item If $j < i$ and $\sigma_i > 0$ then transfer a unit
                          charge from $i$ to $j,$ that is to say set
                          $\sigma_i := \sigma_i-1$ and $\sigma_j := \sigma_j+1.$
                    \item If $j=i$ and $\sigma_i > 0$ then remove a unit charge
                          from $i$ setting $\sigma_i := \sigma_i - 1.$
                    \item If $j > i$ then no update occurs.
                   \end{itemize}
\end{description}
In other words, at each step of the dynamics a charge transfer attempt is made between
two random vertices, which is succesful whenever the source vertex has a higher rank
than the target vertex. Whenever a self-transfer is attempted, a unit charge is removed
from the system (charge leak occurs), although another natural interpretation is that
the evolution of the charge unit terminates at this point and the charge remains stored
forever at the vertex considered rather than being removed from the system, which 
makes {\bf (WTA)} into a charge-conserving dynamics \--- these interpretational issues,
which become important when discussing precise technical relationships between the original
neural network model and its winner-take-all approximation, see Piersa,
Piekniewski \& Schreiber (2010), fall beyond the scope of the present
mathematically oriented article.
 The updates in this dynamics are performed until there are no
more charge units evolving in the system, that is to say $\sigma_i = 0$ for all
$i=1,\ldots,n.$ With each instance of such an evolution
we associate in a natural way its {\it charge-flow network}, also referred to as the
{\it spike-flow network} due to its interpretation in the context of spiking neural
networks as originally considered in Piekniewski \& Schreiber (2008). The charge-flow
network is an undirected graph with multiple edges, where the edge multiplicity
$A^{n,m}_{ij} = A^{n,m}_{ji}$ between $i,j,\; i \geq j,$ is given by the number of
charge units transferred from $i$ to $j$ in the course of the WTA dynamics. Conforming to
the usual terminology, the random symmetric matrix $(A^{n,m}_{ij})_{i,j=1,\ldots,n}$
will be called the adjacency matrix of the charge-flow network in the sequel. Moreover,
the number of charge transfers away from vertex $i,$ that is to say
$\sum_{j \leq i} A^{n,m}_{ij},$ will be called the {\it out-degree} of $i$
and, likewise, the number $\sum_{i \geq j} A^{n,m}_{ij}$ of charge transfers to
vertex $j$ will be called its {\it in-degree} whereas the sum of out- and in-degree
will be called the degree of the vertex. It can be shown, see Theorem 1
in Piekniewski \& Schreiber (2008), whose semi-rigorous proof can easily be
brought to full rigour (which falls beyond the scope of the present work though),
that with overwhelming probability the charge-flow network is asymptotically
scale free with exponent $2$ as $n\to\infty,$  that is to say the in- and out-degrees
of its vertices follow asymptotically a power law with exponent $2,$ see ibidem
for further details.

It is convenient and natural for our further purposes to consider the WTA evolutions
for different values of $n,$ and hence also their corresponding charge flow matrices
$(A^{n,m})_{n \geq 1,\; m = \lfloor \alpha n \rfloor},$ {\it coupled on a common
probability space}, say  $({\Bbb P},\Omega,\Im),$ as follows. For each $n' > n$ the
WTA dynamics on $\{ 1, \ldots, n \}$ is obtained from that on $\{ 1,\ldots, n' \}$
by  
\begin{itemize}
 \item Numbering from $1$ to $m'=\lceil \alpha n'\rceil$ the charge units assigned
         to vertices in $\{ 1,\ldots ,n'\}$ and constructing the restricted initial charge
         assignment for $\{1,\ldots,n \}$ by assigning each among the initial
         $m = \lceil \alpha n \rceil$ units to the first vertex in $\{1,\ldots, n\}$ it hits in
         the course of its extended evolution in $\{ 1,\ldots, n'\}.$ 
\item Letting the evolution of the $m$ charge units in $\{1,\ldots,n\}$ arise as the restriction of
         the corresponding dynamics of the initial $m$ among the $m'$ charge units in $\{1,\ldots,n'\}$ 
         after reaching the set $\{1,\ldots,n\}.$ 
 \end{itemize}
 It is clear that this yields a consistent coupling for all $n \geq 1$ and all our probabilistic
 statements in the sequel shall assume this coupling without a further mention. Note in particular
 that we have almost surely  $(A^{n,m})_{ij} \leq (A^{n',m'})_{ij}$ with $n\leq n', m \leq m'$ and
 $i,j \leq n,$ which allows us to interpret the charge flow graph for $n$ as a subgraph of that for
 $n' \geq n.$ 

The objects in focus of our interest in the present paper are the
(non-normalised!) empirical spectral measures of 
$(A^{n,m}_{ij}),\; m = \lfloor \alpha n \rfloor,$
\begin{equation}\label{SPM}
 \mu_{n,m} := \sum_{i=1}^n \delta_{\lambda_i \slash n}, 
\end{equation}
where $\lambda_1 \geq \lambda_2 \geq \ldots \geq \lambda_n$ are
the eigenvalues of $A^{n,m}$ repeated according to their multiplicities,
note that all $\lambda_i$ are real
numbers because $A^{n,m}$ is self-adjoint. Clearly, the total mass
of $\mu_{n,m}$ is $n,$ but as will be seen in the sequel and as
reflecting the power-law scaling properties of the charge-flow graph,
the random measure $\mu_{n,m}$ with arbitrarily high probability puts
almost all its mass in  neighbourhoods of $0,$ corresponding to the
overwhelming majority of low degree vertices, even though the spectral
radius of $A^{n,m}$ is asymptotically of order $\Theta(n).$ In fact,
we shall show that the mass which $\mu_{n,m}$ puts outside the neighbourhoods
of $0$ is bounded and that, with $n\to\infty$ and 
$m = \lfloor \alpha n \rfloor,$ the random measures $\mu_{n,m}$ converge
almost surely to a non-trivial limit away from $0$ in the sense specified below.

We say that a sequence $\zeta_n$ of Borel measures on ${\Bbb R}$ converges
weakly {\it away from zero} to a Borel measure $\zeta$ on ${\Bbb R}$ iff
$\lim_{n\to\infty} \int f d\zeta_n = \int f d\zeta$ for all bounded
continuous $f : {\Bbb R} \to {\Bbb R}$ which vanish in some neighbourhood 
of zero. To identify the weak limit away from zero for $\mu_{n,m}$
consider the following trace class operator $M: l_2 \to l_2$ on the
space of square-integrable sequences, given by
\begin{equation}\label{MMM}
 [M(a_1,a_2,\ldots)]_i = \sum_{j=1}^{\infty} \frac{a_j}{(i \vee j)^2}.
\end{equation}
Observe that $M$ is symmetric and Hermitian positive as corresponding to
the covariance matrix of $W_{1/i^2},\;i=1,2,\ldots,$ with $W$ standing
for the standard Brownian motion. To get the required trace class property
use that $\sum_{i} 1/i^2 < \infty$ and apply Theorem 2.12 in Simon (2005),
see also ibidem and Section X.3 in Kato (1976) for general theory of trace
class operators.
In particular, the spectrum $\Sigma(M)$ of $M$ is a countable subset
of ${\Bbb R}_+ \cup \{ 0 \}$ with $0$ as its only accumulation point
and each $\lambda \in \Sigma(M),\; \lambda \neq 0,$ is an eigenvalue
of $M.$ Zero belongs to the spectrum as an approximative rather
than proper eigenvalue and, moreover, all eigenvalues of $M$ are
simple. Both these facts are easily checked by writing down the
eigenequation $\lambda a_k = [M(\bar a)]_k$ which yields 
$\lambda(a_{k+1}-a_k) = (1/(k+1)^2 - 1/k^2) \sum_{i=1}^k a_i,\; k \geq 1$ 
\--- clearly the solution to this linear difference equation is unique
up to multiplicative constant for all $\lambda$ and identically zero
for $\lambda=0.$
We set 
\begin{equation}\label{MUINFTY}
 \mu_{\infty} := \sum_{\lambda \in \Sigma(M) \setminus \{ 0 \}} \delta_{\lambda}.
\end{equation}
Our first result states that
\begin{theorem}\label{GLOWNEDYSKR}
 Put $m := \lfloor \alpha n \rfloor.$ Then, with probability one, the sequence
 of random measures $\mu_{n,m}$ converges weakly away from $0$ to $\mu_{\infty} \circ (\alpha)^{-1}$
 as $n\to\infty,$ where $(\alpha)(x) = \alpha x$ stands for the operation of multiplication by $\alpha.$ 
\end{theorem} 
The problem with this theorem, apart from the fact that we are unable to explicitly
determine $\Sigma(M)$ and thus $\mu_{\infty},$ is that it is {\it not robust} with
respect to small modifications of the dynamics, especially for low vertex ranks,
which would have an immediate and non-negligible effect on the operator $M$ and
its spectrum. In particular, the technical issues discussed in the definition of
the {\bf (WTA)} dynamics above and related to the question how to deal with
self-transefers (to regard them charge-leak or charge-freezing events or perhaps
to forbid them at all) do non-trivially impact the limit behaviour of the spectral
measures $\mu_{n,m}.$ 
This is an undesirable situation in our applications to neural nets in the
set-up of Piekniewski \& Schreiber (2008) where the local behaviour of
recurrent neural networks is only approximately driven by the WTA dynamics
and it is at the level of the large-scale global behaviour that we believe
this approximation to yield reliable results. On the other hand, this is
also an unavoidable situation in our present setting, because the spectrum
of the spike-flow graphs is strongly affected by its few highest-degree
vertices.

To get more universal results we need to change somewhat our setting and to
concentrate on medium degree vertices, cutting off those of highest degree
and obtaining theorems characterising the {\it typical architecture} of
the spike-flow graph rather than the individual behaviour of its 
{\it highest order elite} which is highly sensitive to dynamic details.   
To this end, for $\epsilon \in (0,1)$ consider the {\it $\epsilon$-truncated
charge-flow graph} where all connections from and to vertices of rank
between $1$ and $\epsilon n$ are removed (with the downward flow direction
these are the highest degree vertices). The resulting random connectivity
matrix of this graph is denoted by $A^{n,m;\epsilon}.$ We are going to study the
spectral measures
\begin{equation}\label{SPME}
 \kappa_{n,m}^{\epsilon} := \sum_{i=1}^n \delta_{ \epsilon \lambda^{\epsilon}_i}
\end{equation}
where $\lambda^{\epsilon}_1 \geq \lambda^{\epsilon}_2 \geq \ldots$ are the 
eigenvalues of $A^{n,m;\epsilon},$ which are clearly real because
$A^{n,m;\epsilon}$ is symmetric (note that at least 
$\lceil \epsilon_n n \rceil$ among these eigenvalues are $0$ due to the
above cut-off).  As already signalled above, this construction has a very natural
interpretation in terms of large scale neural network modeling purposes in
Piekniewski \& Schreiber (2008) and Piersa, Piekniewski \& Schreiber (2010)
where the effective {\it statistical structure} of the charge flow graph is predominantly
studied at the level of moderate and reasonably high but not highest elite units which
are themselves considered from a somewhat different angle, see e.g. the discussion on
competing {\it basins of attraction} of elite nodes in Section VII.A. of Piersa, Piekniewski
\& Schreiber (2010) for further details. 

To proceed, consider the trace class integral operator
$K:L_2([1,\infty)) \to L_2([1,\infty))$ given by
\begin{equation}\label{KKK}
 [K f](t) = \int_1^{\infty} \frac{f(s)}{(s \vee t)^2} ds.
\end{equation}
As in case of $M$ in (\ref{MMM}) above, also here $M$ is Hermitian positive
as the covariance operator of  $t \mapsto W_{1/t^2},\; t \geq 1,$ and thus
the required trace class property follows by Theorem 2.12 in Simon (2005)
because the trace integral $\int_1^{\infty} 1/t^2 dt$ converges, see
also Example X.1.18 in Kato (1976). 
In particular, the spectrum of $K$ consists of a countable set
of isolated positive eigenvalues accumulating at $0.$ Zero belongs
to the spectrum as an approximative rather than proper eigenvalue.
In contrast to $M$ here we are able to explicitly determine the
spectrum of $K$ though.
\begin{lemma}\label{SPECK}
 All eigenvalues of $K$ are simple and strictly positive. Moreover, for $\lambda > 0$ we have 
 $$\lambda \in \Sigma(K) \Leftrightarrow J_1\left(\frac{2\sqrt{2}}{\sqrt{\lambda}}\right)
 = 0$$
 where $J_1(\cdot)$ is the Bessel J-function of order $1.$
\end{lemma}  
We put
\begin{equation}\label{KAPINFTY}
 \kappa_{\infty} := \sum_{\lambda \in \Sigma(K) \setminus \{ 0 \}} \delta_{\lambda}.
\end{equation}
Choose a sequence 
$(\epsilon_n)_{n=1}^{\infty},$ in the sequel often required to satisfy
\begin{equation}\label{EPSWAR}
 \lim_{n\to\infty} n \epsilon_n = +\infty \;\; \mbox{ and there exists $\delta > 0$
  such that } \; \lim_{n\to\infty} n^{1+\delta} \epsilon^2_n = 0. 
\end{equation}
Our second main result is 
\begin{theorem}\label{GLOWNECGL}
 Put $m := \lfloor \alpha n \rfloor$ and let $\epsilon_n$ be as in (\ref{EPSWAR}).
 Then, with probability one, the sequence of random measures $\kappa_{n,m}^{\epsilon_n}$
 converges weakly away from $0$ to $\kappa_{\infty} \circ (\alpha)^{-1}$ as $n\to\infty.$
\end{theorem} 
 The interpretation of the first condition in (\ref{EPSWAR}) is rather clear in this
 context \--- we want the cut-off rank $\epsilon_n n$ to move towards $+\infty$ as
 $n$ does. The second condition in (\ref{EPSWAR}) is perhaps less intuitive and its
 origin will be explained in the discussion following the proof of Theorem \ref{GLOWNECGL}.

 Upon inspecting its proof, Theorem \ref{GLOWNECGL} is easily seen to be insensitive
 to local dynamic modifications, such as these discussed following the formulation
 of Theorem \ref{GLOWNEDYSKR}, whose impact is only sensed by eigenvalues in close
 neighbourhoods of $0.$ This is an important good news from the viewpoint
 of our envisioned applications to large scale neural networks. 

 We conclude this section by one further important remark. It is known, see (9.57)
 in Temme (1996), that $k$-th zero of the Bessel function $J_1$ is asymptotic to
 $1/4 + k \pi$ as $k \to \infty.$ Consequently, by Theorem \ref{GLOWNECGL},
 the $k$-th eigenvalue of $\kappa_{n,m}^{\epsilon_n}$ asymptotically 
 approaches $\frac{8 \alpha}{\pi^2 k^2}$ for large $k.$ This means that 
 the spectral measures $\kappa_{n,m}^{\epsilon_n}$ asymptotically reproduce
 the power law with exponent $2$ as gouverning the degree distribution of
 the charge flow graph, see Piekniewski \& Schreiber (2008). This is rather
 natural since the large eigenvalues of the considered adjacency graph
 are due to its large degree vertices. 

\section{Proofs}

\subsection{Proof of Theorem \ref{GLOWNEDYSKR}}
 The proof of our Theorem \ref{GLOWNEDYSKR} uses the convergence of moments of spectral
 measures $\mu_{n,m}$ which admit convenient representation as the traces of respective
 powers of the adjacency matrix of the considered charge flow graph. We put
 \begin{equation}\label{KMOM}
  M_{k,n} := \int_{\Bbb R} \lambda^k d\mu_{n,m}(\lambda).
 \end{equation}
 First we shall show that the desired convergence of moment expectations holds:
 \begin{lemma}\label{TW1}
  With the notation above we have for $k \geq 1$
  $$
  \lim_{n \to \infty} {\Bbb E}M_{k,n} = \alpha^k \int \lambda^k \mu_{\infty}(d\lambda). 
  $$ 
 \end{lemma}
 Next, applying appropriate measure concentration techniques, we will use
 Lemma \ref{TW1} to show that
 \begin{corollary}\label{TW1AS}
  We have almost surely
  $$
   \lim_{n \to \infty} M_{k,n} = \alpha^k \int \lambda^k \mu_{\infty}(d\lambda).
  $$  
  \end{corollary}
  Finally, applying Corollary \ref{TW1AS} we will complete the proof of
  Theorem \ref{GLOWNEDYSKR} by standard argument.

  \paragraph{Proof of Lemma \ref{TW1}}
 To calculate ${\Bbb E}M_{k,n}$ we write first
 \begin{equation}\label{TRR}
  {\Bbb E}M_{k,n} = {\Bbb E}\tr([A^{n,m}]^k)/n^{k}.
 \end{equation}
 As already indicated in the construction of our standard coupling between the WTA dynamics for different 
 system sizes,  we adopt the convenient convention of numbering from $1$ to $m$ the charge units 
 present in the system.
 Under this convention, whenever a transfer is made from vertex $i$ to $j,$ the
 number of unit to be transferred is chosen in some deterministic way among
 the numbers ascribed to units stored at $i,$ for instance the lowest/highest
 or the first/last arrived one. Consequently, recalling the dynamics of the
 system we get from (\ref{TRR})
 \begin{equation}\label{UREPR}
  n^k {\Bbb E}M_{k,n} = \sum_{l_1=1}^m \ldots \sum_{l_k=1}^m \sum_{U_1=1}^n \ldots \sum_{U_k=1}^n
  {\Bbb P}\left( {\cal T}(U_1,U_2;l_1) \cap {\cal T}(U_2,U_3;l_2) \ldots \cap
  {\cal T}(U_k,U_1;l_k) \right), 
 \end{equation}
 where ${\cal T}(U_{i},U_{i+1};l_i)$ stands for the event that the $l_i$-th charge unit
 was directly transferred between vertices $U_i$ and $U_{i+1},$ either from $U_i$ to
 $U_{i+1}$ or in the opposite direction, in the course of the system evolution.
 To proceed, we split the RHS of (\ref{UREPR}) into a sum of two terms:
 \begin{itemize}
  \item $S_k$ given as the sum of the RHS terms of (\ref{UREPR}) for which all $l_i$'s
        are different,
  \item $R_k$ given as the sum of the remaining terms in the RHS of (\ref{UREPR}), that is
        to say these where at least two $l_i$'s coincide.   
 \end{itemize}
 We evaluate $S_k$ first, and then we show that $R_k$ is of a smaller order and thus
 asymptotically negligible. Since the sequences of vertices visited by different
 charge units on their way to $1$ are independent, we have
 $$ S_k = \sum_{\stackrel{l_i \in \{1,\ldots,m\},\; i=1,\ldots,k}{\forall_{i \neq j} l_i \neq l_j}}
    \sum_{U_1=1}^n \ldots \sum_{U_k=1}^n  {\Bbb P}({\cal T}(U_1,U_2;l_1))
    {\Bbb P}({\cal T}(U_2,U_3;l_2)) \ldots {\Bbb P}({\cal T}(U_k,U_1;l_k)) =
 $$
 \begin{equation}\label{SKW1}
    m (m-1) \ldots (m-k+1) \sum_{U_1=1}^n \ldots \sum_{U_k=1}^n  {\Bbb P}({\cal T}(U_1,U_2;l_1))
    {\Bbb P}({\cal T}(U_2,U_3;l_1)) \ldots {\Bbb P}({\cal T}(U_k,U_1;1)),
 \end{equation}
 with the last equality due to the fact that the evolutions of all charge units 
 coincide in law as following the same dynamic rules. To evaluate the probability
 of ${\cal T}(U_i,U_{i+1};1)$ assume with no loss of generality that $U_{i+1} \leq U_i.$
 Then, since the number of the next vertex to be visited by a unit charge in the course
 of its WTA evolution is uniform among the numbers not exceeding the current vertex number,  
 we have
 \begin{equation}\label{DZIEL1} 
  {\Bbb P}({\cal T}(U_i,U_{i+1};1)) = \frac{1}{U_i} {\Bbb P}({\cal T}(U_i;1)),
 \end{equation}
 where ${\cal T}(U_i;1)$ is the event that $1$-st charge unit has
 visited the vertex $U_i$ on its way towards $1.$ Now, to find 
 ${\Bbb P}({\cal T}(U_1;1))$ note that, by standard extreme value
 theory for record statistics as discussed e.g. in Subsection 4.1
 in Resnick (1987), the sequence of different vertices $V_1 > V_2,\ldots$ visited
 by a charge unit coincides in law with the sequence
 \begin{equation}\label{POISREPR}
  \lceil n \exp(-\eta_1) \rceil,\lceil n \exp(-\eta_2) \rceil, \ldots,
 \end{equation}
 where $\eta_i$ is the $i$-th consecutive point of a
 homogeneous Poisson point process of intensity $1$ on ${\Bbb R}_+$ conditioned
 on not having more than one point in any of the intervals 
 $[-\log(U/n),-\log((U-1)/n)),\; U \in \{1,\ldots,n\}$ under the convention
 that $\log 0 = -\infty.$ Consequently, ${\Bbb P}({\cal T}(U_1;1))$ coincides
 with the probability that some Poisson point $\eta_i$ falls into
 $[-\log(U_1/n),-\log((U_1-1)/n))$ which is $1-\exp(-\log(U_1/n)+\log((U_1-1)/n)) =
 1 - \frac{U_1-1}{U_1} = 1 / U_1.$ Thus, we conclude from (\ref{DZIEL1}) that 
 \begin{equation}\label{TUU}
  {\Bbb P}({\cal T}(U_i,U_{i+1};1)) = \frac{1}{(U_i \vee U_{i+1})^2}
 \end{equation}
 and hence, by (\ref{SKW1}),
 \begin{equation}\label{CONCL}
  S_k = m (m-1) \ldots (m-k+1) \sum_{U_1=1}^n \ldots \sum_{U_k=1}^n  
  \prod_{i=1}^k \frac{1}{(U_i \vee U_{i+1})^2}
 \end{equation}
 with the convention that $U_{k+1} = U_1.$
 Further, we want to estimate the contribution brought by the extra term
 $R_k.$ We claim that
 \begin{equation}\label{CONCL2}
   S_k \leq S_k + R_k \leq k! \left(^m_k\right) \sum_{U_1=1}^n \ldots \sum_{U_k=1}^n  
  \prod_{i=1}^k \frac{1}{U_i \vee U_{i+1}} \left[\frac{1}{U_i \vee U_{i+1}} + 1/\Theta(m) \right].
 \end{equation}
 Indeed, whenever $l_{i+1} = l_i,$ the events ${\cal T}(U_i,U_{i+1};l_i)$ and
 ${\cal T}(U_{i+1},U_{i+2};l_{i+1})$ are no more independent and in fact can
 only co-occur if $U_{i+1}$ lies between $U_i$ and $U_{i+2},$ i.e. $U_i \leq U_{i+1}
 \leq U_{i+2}$ or $U_i \geq U_{i+1} \geq U_{i+2},$ for otherwise one transfer would
 have two different sources or two different destinations. Thus, if we proceeded
 as in our derivation of (\ref{SKW1}) for $S_k,$ we would lose the factor
 $\frac{1}{U_{i+1}}$ corresponding to ${\Bbb P}({\cal T}(U_{i+1};l_{i+1}))$
 since ${\cal T}(U_{i+1};l_{i+1}) = {\cal T}(U_{i+1};l_i).$ We would get
 $1$ instead, but on the other hand we would lose the summation over $l_{i+1},$
 which is now $l_i.$ This means losing one of the $k$ prefactors of order $\Theta(m)$
 as present in the RHS of (\ref{CONCL}) above or, equivalently, keeping summation
 over a dummy variable $l_{i+1}'$ not to lose any prefactors, but with the lost
 factor $\frac{1}{U_{i+1}}$ replaced by $1/\Theta(m)$ for each instance of $l_{i+1}'.$
 This justifies  (\ref{CONCL2}) as required. 
 Thus, recalling that $m = \lfloor \alpha n \rfloor,$
 using (\ref{UREPR}) and combining (\ref{CONCL}) and (\ref{CONCL2})
 we obtain
 \begin{equation}\label{WYNIK}
  \lim_{n\to\infty} {\Bbb E} M_{n,k} = \alpha^k \sum_{U_1=1}^{\infty} \ldots \sum_{U_k=1}^{\infty}  
  \prod_{i=1}^k \frac{1}{(U_i \vee U_{i+1})^2} 
 \end{equation}
 with the convergence of the RHS series easily verified. Finally, recalling (\ref{MMM}),
 using (\ref{WYNIK}) and the trace class properties of $M^k$ yields
 $$ \lim_{n\to\infty} {\Bbb E} M_{n,k} = \alpha^k \tr M^k $$ 
 which completes the proof of Lemma \ref{TW1} in view of the spectral measure definition
 (\ref{MUINFTY}). $\Box$

\paragraph{Proof of Corollary \ref{TW1AS}}
 We begin by considering a modified version of our basic WTA dynamics,
 which is better suited for an application of measure concentration results
 whereas with overwhelming probability its resulting charge-flow graph does
 coincide with the original winner-take-all network. The modification
 is that whenever on its way towards $1$ a charge unit makes more than
 $n^{1/3}$ jumps, then it is forced to make its final jump directly to $1$
 rather than further following the usual dynamics. By our Poisson representation
 (\ref{POISREPR}) of single charge unit evolution the number of jumps
 made on the way to $1$ behaves asymptotically as mean $\log n$ Poisson
 random variable $\Po(\log n).$ Consequently, the probability that the
 number of jumps of an individual charge unit exceeds $n^{1/3}$ is
 not larger than $\exp\left(-\frac{n^{1/3}}{4}\log(n^{1/3}/2)\right),$
 see e.g. Shorack \& Wellner (1986), p. 485. Thus, since
 the overall number of charge units is $m = \lfloor \alpha n \rfloor,$
 the probability that {\it any} individual charge unit makes more than $n^{1/3}$
 jumps is still of order $\exp(-\Theta(n^{1/3} \log n)).$
 Writing $\hat{A}^{n,m}$ for the adjacency matrix under
 the modified dynamics we have therefore
 \begin{equation}\label{AHATA}
  {\Bbb P}(\hat{A}^{n,m} \neq A^{n,m}) \leq \exp(-\Theta(n^{1/3}\log n)).
 \end{equation}

 To complete the proof we shall proceed by induction in $k.$ Assume first
 that $k=1$ and note that $\tr(\hat{A}^{n,m})$ is a $1$-Lipschitz function of 
 $\hat{A}^{n,m}$ under the $l_1$-norm on ${\Bbb R}^{n \times n}.$
 Consider now the operation of replacing the evolution of a single charge
 unit under the modified dynamics by some other evolution with at most $n^{1/3}$ jumps.
 Let $B$ be the difference matrix between the new and the original adjacency
 matrices $\hat{A}^{n,m}.$ Clearly, $B$ has at most $4n^{1/3}$ non-zero entries,
 all of which are ones or minus ones. Thus, such an operation may change 
 $\tr(\hat{A}^{n,m})$ by at most $4n^{1/3}$ and, consequently, 
 $\tr(\hat{A}^{n,m}/n)$ by at most $4n^{-2/3}.$  
 Recalling that $\hat{A}^{n,m}$ is a function of the evolutions of $m$ individual
 charge units wich are independent, and using standard measure concentration results
 for Lipschitz functions of independent entries, see Corollary 1.17 in Ledoux (2001),
 we conclude that
 \begin{equation}\label{DLAK1}
  {\Bbb P}(|\tr(\hat{A}^{n,m}/n) - {\Bbb E}\tr(\hat{A}^{n,m}/n)| \geq t)
  \leq 2 \exp\left(-\frac{t^2}{\Theta(m n^{-4/3})} \right) = \exp(-\Theta(t^2 n^{1/3}))
 \end{equation}
 because $m = \lfloor \alpha n \rfloor.$
 With $t := 1/(\log n)$ relation (\ref{DLAK1}) becomes 
 \begin{equation}\label{DLAK1P}
  {\Bbb P}(|\tr(\hat{A}^{n,m}/n) - {\Bbb E}\tr(\hat{A}^{n,m}/n)| \geq 1/(\log n))
  \leq \exp(-\Theta(n^{1/3} (\log n)^{-2})). 
 \end{equation}
 Combining (\ref{DLAK1P}) with (\ref{AHATA}) above yields now
 $$ {\Bbb P}(|\tr(A^{n,m}/n) - {\Bbb E}\tr(A^{n,m}/n)| \geq 1/(\log n))
    \leq \exp(-\Theta(n^{1/3} (\log n)^{-2})) $$
 whence the assertion of the corollary for $k=1$ trivially follows by the Borel-Cantelli lemma.

 To proceed with our inductive argument, for technical convenience we slightly extend our assertion
 for $k \geq 2$ and we show that both 
 \begin{equation}\label{HIPIND}
  {\Bbb P}(|\tr([\hat{A}^{n,m}/n]^k) - {\Bbb E}\tr([\hat{A}^{n,m}/n]^k)| \geq 1/(\log n))
  \leq \exp(-\Theta(n^{1/3} (\log n)^{-2}))
 \end{equation}
  and
 \begin{equation}\label{HIPINDabs}
    {\Bbb P}(|\tr([{\rm abs}(\hat{A}^{n,m})/n]^k) - {\Bbb E}\tr([{\rm abs}(\hat{A}^{n,m})/n]^k)| \geq 1/(\log n))
  \leq \exp(-\Theta(n^{1/3} (\log n)^{-2}))
 \end{equation}
 hold for all $k \geq 2, $ with the absolute value matrix ${\rm abs}(\hat{A}^{n,m})$ understood here
 in the usual spectral sense (the same eigenvectors, eigenvalues replaced by absolute values). Assuming
 that (\ref{HIPIND}) and (\ref{HIPINDabs}) have already been established for $k-1$ (unless $k=2$
 where we only assume (\ref{HIPIND}) to hold) we define an 
 auxiliary modified trace functional $\hat{\tr}_k(\cdot),\; k \geq 2,$ by putting for an $n \times n$ matrix $A$
 \begin{enumerate}
  \item If $k=2$ and 
    \begin{equation}\label{RNIEAAdla2}
     \tr(A) \leq 2 \alpha^{k-1} \int \lambda \mu_{\infty}(d\lambda)
    \end{equation}
    then $\hat{\tr}_k(A) := \tr(A^k),$
  \item If $k\geq 3$ and $k$ is odd and 
   \begin{equation}\label{RNIEAAdlaParz}
    \tr(A^{k-1}) \leq 2 \alpha^{k-1} \int \lambda^{k-1} \mu_{\infty}(d\lambda)
    \end{equation}
    then $\hat{\tr}_k(A) := \tr(A^k),$
  \item If $k\geq 3$ and $k$ is even and
   \begin{equation}\label{RNIEAAdlaNieparz}
   \tr({\rm abs}(A)^{k-1}) \leq 2 \alpha^{k-1} \int \lambda^{k-2} + \lambda^k \mu_{\infty}(d\lambda)
  \end{equation}
  then $\hat{\tr}_k(A) := \tr(A^k),$
 \item Otherwise, define $\hat{\tr}_k(A) := \tr(\tilde{A}^k)$ where $\tilde{A}$ is
         the metric projection of $A$ onto the set 
         ${\cal A}_{\mu_{\infty}} = {\cal A}[k,n,\alpha,\mu_{\infty}]$ given as
         \begin{enumerate}
          \item the set of $n\times n$ matrices satisfying (\ref{RNIEAAdla2}) if $k=2$  (as in case 1.)
          \item the set of $n\times n$ matrices satisfying (\ref{RNIEAAdlaParz}) if $k \geq 3$ and $k$ is odd
           (as in case 2.)
          \item the set of $n \times n$ matrices satisfying (\ref{RNIEAAdlaNieparz}) if $k \geq 3$  and $k$ is even
           (as in case 3.)
         \end{enumerate}
         Note that by Klein's lemma, see e.g. Lemma 6.4 in Guionnet (2009), the set 
         ${\cal A}_{\mu_{\infty}}$ is convex and closed and the matrix $\tilde{A}$ is simply the
         matrix in ${\cal A}_{\mu_{\infty}}$ minimising the Euclidean distance to
         $A$ in ${\Bbb R}^{n \times n}.$
 \end{enumerate}
 This somewhat technical definition has a very simple interpretation: the modified trace
 functional $\hat{\tr}_k(\cdot)$ coincides with the usual trace of $A^k$
 provided that the corresponding trace of ${\rm abs}(A)^{k-1}$ is not too large,
 otherwise the modified trace is defined
 as the trace of $\tilde{A}^k$ where $\tilde{A}$ is a version of the matrix $A$
 projected onto an appropriate convex set ${\cal A}_{\mu_{\infty}}$ so that the
 trace of its $(k-1)$-th power does not exceed the corresponding controllable threshold
 given by the RHS of (\ref{RNIEAAdla2},\ref{RNIEAAdlaParz}) and (\ref{RNIEAAdlaNieparz})
 respectively, we denote this threshold by $\tau[k,\alpha,\mu_{\infty}]$ for reference below.
 The extra auxiliary relation (\ref{HIPINDabs}) is needed to ensure the convexity of
 ${\cal A}_{\mu_{\infty}}$ for $k$ even.

The further argument is quite standard now: the above modified trace functional coincides with the
 original one with overwhelming probability and at the same time it is well behaved as admitting
 well controllable oscillations and thus is suitable for usual measure concentration techniques. 
 Indeed, using (\ref{AHATA}) and applying Lemma \ref{TW1} combined with the observation that
 $|\lambda|^{k-1} \leq \lambda^{k-2} + \lambda^k$ for $k \geq 2$ even, we conclude from (\ref{HIPIND})
 and (\ref{HIPINDabs}) for $k-1$ that 
 \begin{equation}\label{MALAROZNICA}
  {\Bbb P}(\tr([\hat{A}^{n,m}/n]^k) \neq \hat{\tr}_k(\hat{A}^{n,m}/n)) \leq
  \exp(-\Theta(n^{1/3} (\log n)^{-2})).
 \end{equation}
 Recalling now that the derivative of $A \mapsto \tr(A^k)$ in the direction of a matrix
 $B$ is given by $k \tr(A^{k-1} B)$ (see e.g. Lemma 6.1 in Guionnet (2009)),
 taking $A := \hat{A}^{n,m}$ and letting $B$ be 
 the evolution replacement difference matrix as above, we conclude by convexity of
 the projection set ${\cal A}_{\mu_{\infty}}$ and upon recalling that $B$ has
 at most $4 n^{1/3}$ non-zero entries, all plus or minus ones, that
 $$ |\hat{\tr}_k([\hat{A}^{n,m}+B]/n) 
   - \hat{\tr}_k([\hat{A}^{n,m}/n])| \leq k n^{-1} 4 n^{1/3} \tau[k,\alpha,\mu_{\infty}] = \Theta(n^{-2/3}). $$ 
 Thus, using again that $\hat{A}^{n,m}$ is a function of the evolutions of $m$ individual
 charge units wich are independent, and applying one more time Corollary 1.17 in Ledoux (2001),
 we obtain
 \begin{equation}\label{DLAKOG}
  {\Bbb P}(|\hat{\tr}_k(\hat{A}^{n,m}/n) - {\Bbb E}\hat{\tr}_k(\hat{A}^{n,m}/n)| \geq 
   1/(\log n)) \leq \exp(-\Theta(n^{1/3} (\log n)^{-2}))
 \end{equation}
 in full analogy to (\ref{DLAK1P}). When combined with (\ref{MALAROZNICA}) this yields
 the required relation (\ref{HIPIND}) for $k.$ The second inductive relation (\ref{HIPINDabs})
 follows in full analogy by using the fact that the derivative of $A \mapsto \tr({\rm abs}(A^k))$
 in the direction of a matrix $B$ is $k \tr({\rm abs}(A)^{k-1} B)$ for $k \geq 2,$ see again e.g.
 Lemma 6.1 in Guionnet (2009). This completes the inductive argument and shows that both
 (\ref{HIPIND}) and (\ref{HIPINDabs}) hold for all $k \geq 2.$ 
 
 Finally, putting (\ref{HIPIND})  together with (\ref{AHATA}) we come to
 \begin{equation}\label{KONCK}
   {\Bbb P}(|\tr([A^{n,m}/n]^k) - {\Bbb E}\tr([A^{n,m}/n]^k)| \geq 1/(\log n))
   \leq \exp(-\Theta(n^{1/3} (\log n)^{-2}))
 \end{equation}
 for all $k \geq 1,$ which completes the proof of Corollary \ref{TW1AS} by a
 straightforward application of the Borel-Cantelli lemma.
 $\Box$

 \paragraph{Completing the proof of Theorem \ref{GLOWNEDYSKR}}

 Having established Corollary \ref{TW1AS} we readily complete the proof
 of Theorem \ref{GLOWNEDYSKR} using that the trace class operator $M$
 has in particular a finite spectral radius
 and resorting to the standard method of moments and classical Carleman's
 criterion, see e.g. Shohat \& Tamarkin (1943), p. 19, applied for the
 measures $\mu'_{n,m}(d\lambda) := \lambda \mu_{n,m}(d\lambda)$ whose
 sequence of moments coincides with that of $\mu_{n,m}$ shifted by one
 \--- this way we conclude that a.s. $\mu'_{n,m}$ converges weakly to
 $\mu'_{\infty}$ with $\mu'_{\infty}(d\lambda) = 
 \lambda \mu_{\infty}(d\lambda)$ whence the desired a.s. weak convergence
 of $\mu_{n,m}$ to $\mu_{\infty}$ away from zero follows.
 $\Box$

 \paragraph{Remarks}
  An intuitive explanation of Theorem \ref{GLOWNEDYSKR} can be provided by noting
  that, in view of (\ref{TUU}), we have
  $$ {\Bbb E} A^{n,m}_{ij} = \frac{m}{(i \vee j)^2} $$
  and the fluctuations of $A^{n,m}_{ij}$ are easily controllable as coming
  from independent evolutions of $m$ charge units. Consequently, $A^{n,m}/n$
  a.s. converges to $\alpha M$ in many reasonably strong senses provided by
  the operator theory
  and thus $\mu_{\infty} \circ (\alpha)^{-1}$ is a natural candidate for
  the limit of spectral measures $\mu_{n,m}.$ This could be a starting point for an alternative
  proof of Theorem \ref{GLOWNEDYSKR}, but presumably much more complicated
  than ours as requiring the use of measure concentration tools in Banach space
  of linear operators endowed with the trace class norm, and then quite involved
  and technical additional considerations relating the convergence of operators
  to spectral measure convergence. In this context, we strongly prefer the
  method of moments as letting us  avoid unnecessary technicalities.

\subsection{Proof of Theorem \ref{GLOWNECGL}}
 As in the proof of Theorem \ref{GLOWNEDYSKR} also here we use the
 convergence of moments. With $m = \lfloor \alpha n \rfloor$ we put
 \begin{equation}\label{KKMOM}
  M^{\epsilon}_{k,n} := \int_{\Bbb R} \lambda^k d\kappa^{\epsilon}_{n,m}(\lambda).
 \end{equation}
 We shall establish the following covergence of expectations first.
 \begin{lemma}\label{KTW1}
  With $\epsilon_n$ such that
  $\lim_{n\to\infty} \epsilon_n n = +\infty$ and 
  $\lim_{n\to\infty} \epsilon_n = 0$  
  we have for $k \geq 1$
  $$
   \lim_{n \to \infty} {\Bbb E} M^{\epsilon_n}_{k,n} = 
   \alpha^k \int \lambda^k \kappa_{\infty}(d\lambda). 
  $$ 
 \end{lemma}
 In analogy to the corresponding step in the proof of Theorem \ref{GLOWNEDYSKR},
 also here the convergence of expectations will be strengthened to a.s. convergence
 using measure concentration. 
 \begin{corollary}\label{KTW1AS}
  Assume that the sequence $\epsilon_n$ satisfies (\ref{EPSWAR}).
  Then we have almost surely
  $$
   \lim_{n \to \infty} M^{\epsilon_n}_{k,n} = \alpha^k \int \lambda^k \mu_{\infty}(d\lambda).
  $$  
  \end{corollary}
  Note that, unlike in Lemma \ref{KTW1}, in Corollary \ref{KTW1AS} we do require 
  the full strength of (\ref{EPSWAR}).
  This corollary will lead us to the desired assertion of Theorem \ref{GLOWNECGL}
  by a standard argument. 
  \paragraph{Proof of Lemma \ref{KTW1}}
   To calculate ${\Bbb E}M^{\epsilon_n}_{k,n}$ write
   $$
   {\Bbb E}M^{\epsilon_n}_{k,n} = \epsilon_n^{k} {\Bbb E}\tr([A^{n,m;\epsilon_n}]^k).
   $$ 
   In full analogy with the corresponding argument leading to (\ref{CONCL}) and
   (\ref{CONCL2}) in the proof of Theorem \ref{GLOWNEDYSKR} above, we obtain
   \begin{equation}\label{NIERCONCL}
     k! \left(^m_k\right) \epsilon_n^{k} 
    \sum_{U_1=\lceil \epsilon_n n \rceil}^n \ldots \sum_{U_k=\lceil \epsilon_n n \rceil}^n
     \prod_{i=1}^k \frac{1}{(U_i \vee U_{i+1})^2}
    \leq {\Bbb E} M^{\epsilon_n}_{k,n} \leq
   \end{equation}
    $$ 
     k! \left(^m_k\right) \epsilon_n^{k}
     \sum_{U_1=\lceil \epsilon_n n \rceil}^n \ldots \sum_{U_k=\lceil \epsilon_n n \rceil}^n  
     \prod_{i=1}^k \frac{1}{U_i \vee U_{i+1}} \left[\frac{1}{U_i \vee U_{i+1}} + 1/\Theta(m)
     \right].       
    $$
    Consequently, as $n\to\infty,$ we have in view of (\ref{NIERCONCL}) 
    $$ {\Bbb E} M^{\epsilon_n}_{k,n} = (1+o(1)) \frac{k! \left(^m_k\right)}{n^k}
       \frac{1}{(\epsilon_n n)^k}
       \sum_{U_1=\lceil \epsilon_n n \rceil}^n \ldots \sum_{U_k=\lceil \epsilon_n n \rceil}^n
       \prod_{i=1}^k \left(\frac{\epsilon_n n}{U_i \vee U_{i+1}}\right)^2. $$
    Substituting $u_i := U_i / (\epsilon_n n),$ recognising appropriate integral sums
    in the RHS and recalling that $m = \lfloor \alpha n \rfloor$ we get therefore
    by our assumptions on $\epsilon_n$
    $$ \lim_{n\to\infty} {\Bbb E} M^{\epsilon_n}_{k,n} = \alpha^k 
       \int_1^{\infty} \ldots \int_1^{\infty} \prod_{i=1}^k \frac{1}{(u_i \vee u_{i+1})^2}
       du_1 \ldots du_k. $$
    Recalling the definition (\ref{KKK}) of $K$ and the trace class properties of $K^k$ this yields
    $$ \lim_{n\to\infty} {\Bbb E} M^{\epsilon_n}_{k,n} = \alpha^k \tr K^k. $$
    This completes the proof of Lemma \ref{KTW1} in view of the spectral measure
    definition (\ref{KAPINFTY}).
    $\Box$ 
  
 \paragraph{Proof of Corollary \ref{KTW1AS}}
  Our argument here goes very much along the same lines as the proof of Corollary \ref{TW1AS}.
  Note first that Lemma \ref{KTW1} is applicable under the assumptions of Corollary
  \ref{KTW1AS} because (\ref{EPSWAR}) does in particular imply the conditions on
  $\epsilon_n$ imposed in the statement of the lemma. 
  Again, we consider a modified version of the WTA dynamics,  the modification
  being that whenever on its way towards $1$ a charge unit makes more than
  $n^{\delta/3}$ jumps, then it is forced to make its final jump directly to $1$
  rather than further following the usual dynamics. Recall that $\delta$ is determined
  by (\ref{EPSWAR}) as assumed in the statement of the corollary. Writing again
  $\hat{A}^{n,m}$ for the adjacency matrix under the modified dynamics we have
  in full analogy with (\ref{AHATA})
  \begin{equation}\label{AHATA2}
   {\Bbb P}(\hat{A}^{n,m} \neq A^{n,m}) \leq \exp(-\Theta(n^{\delta/3}\log n)).
  \end{equation}
  In analogy to the proof of Corollary \ref{TW1AS}, also here we consider
  the operation of replacing the evolution of a single charge  unit under
  the modified dynamics by some other evolution with at most $n^{\delta/3}$ jumps.
  Denoting by $B$ be the difference matrix between the new and the original
  adjacency matrices $\hat{A}^{n,m}$ we see that $B$ has at most $4n^{\delta/3}$
  non-zero entries, all of which are ones or minus ones. This observation
  puts us again in a position to apply measure concentration results for
  Lipschitz functionals with respect to product measures, nearly verbatim
  following the respective lines of the inductive argument for Corollary 
  \ref{TW1AS}. Note that in our present set-up the modified trace functional
  $\hat{\tr}_k$ involves projections onto the convex set ${\cal A}_{\kappa_{\infty}}$
  defined in full analogy to the corresponding ${\cal A}_{\mu_{\infty}}.$ 
  Moreover, $\hat{A}^{n,m}/n$ in the proof of Corollary \ref{TW1AS}
  is replaced by $\epsilon_n \hat{A}^{n,m}$ here due to the different scaling.
  This way, in analogy to (\ref{HIPIND}), we conclude that
  \begin{equation}\label{HIPIND2}
  {\Bbb P}(|\tr([\epsilon_n \hat{A}^{n,m}]^k) - {\Bbb E}\tr([\epsilon_n \hat{A}^{n,m}]^k)| \geq 1/(\log n))
  \leq \exp(-\Theta(n^{\delta/3} (\log n)^{-2}))
  \end{equation}
  for all $k \geq 1.$
  Using (\ref{AHATA2}) we get 
  $$
  {\Bbb P}(|\tr([\epsilon_n A^{n,m}]^k) - {\Bbb E}\tr([\epsilon_n A^{n,m}]^k)| \geq 1/(\log n))
  \leq \exp(-\Theta(n^{\delta/3} (\log n)^{-2}))
  $$
  in analogy to (\ref{KONCK}), whence the assertion Corollary \ref{KTW1AS} follows by the
  Borel-Cantelli lemma.
  $\Box$

  \paragraph{Completing the proof of Theorem \ref{GLOWNECGL}}
 
  Since the trace class operator $K$ has in particular a finite spectral radius,
  the desired assertion of Theorem \ref{GLOWNECGL} follows now readily in view of
  Corollary \ref{KTW1AS} by the standard method of moments and Carleman's criterion,
  see e.g. p. 19 in Shohat \& Tamarkin (1943), used in analogy to the corresponding
  proof-completing paragraph for Theorem \ref{GLOWNEDYSKR}. 
 $\Box$
 
 \paragraph{Justification of condition (\ref{EPSWAR})}
  We note at this point that, intuitively speaking, the independent contributions to the
  random matrix $A^{n,m}$ brought by each of the $m$ evolving charge units should bring
  respective variance contributions to the trace $\tr([\epsilon_n A^{n,m}]^k)$ of the
  order $\epsilon_n^2 \log n$ per unit ($\log n$ is the order of number of unit charge
  jumps before leaking out from the system), which sums up to order 
  $\Theta(m \log n \epsilon_n^2) = \Theta(n \log n \epsilon_n^2)$
  upon taking all units into account. Therefore it is natural to require that
  $n \epsilon_n^2$ converges to $0$ faster than $1/(\log n),$ which is roughly the
  content of the second condition in (\ref{EPSWAR}), for otherwise we should not
  hope for a  deterministic limit of $\tr([\epsilon_n A^{n,m}]^k)$ 
  and thus of $\mu_{n,m;\epsilon_n}$ as $n\to\infty.$
  This informal observation should be regarded as a justification for
  (\ref{EPSWAR}) rather than as a mathematical statement though.

\subsection{Proof of Lemma \ref{SPECK}}
 Assume that $\phi \in L_2([1,\infty))$ is a non-zero eigenfunction of the operator
 $K$ corresponding to some eigenvalue $\lambda \geq 0.$ The corresponding eigenequation
 reads
 \begin{equation}\label{RNIEWL}
  \lambda \phi(t) = \frac{1}{t^2} \int_1^t \phi(s) ds + \int_t^{\infty} \frac{1}{s^2} \phi(s) ds.
 \end{equation}
 Since the RHS is an application of the integral operator with a well-behaved kernel,
 both sides are readily seen to be differentiable and the differentiation yields
 $$ \lambda \phi'(t) = \frac{1}{t^2} \phi(t) - \frac{1}{t^2} \phi(t) - \frac{2}{t^3} \int_1^t \phi(s) ds
    = - \frac{2}{t^3} \int_1^t \phi(s) ds. $$
 Putting 
 $$ \Psi(t) := \int_1^t \phi(s) ds $$
 we get the differential equation
 \begin{equation}\label{RNIER}
  \lambda \Psi''(t) = - \frac{2}{t^3} \Psi(t)
 \end{equation}
 with the initial condition
 \begin{equation}\label{WARPOCZ}
  \Psi(1) = 0.
 \end{equation}
 The solution to this equation is $\Psi \equiv 0$ for $\lambda = 0$ which shows that
 $0$ is not an eigenvalue and, for $\lambda \neq 0,$
 \begin{equation}\label{PSIWYN}
  \Psi(t) = \sqrt{t} \left( C_1 J_1\left( \frac{2\sqrt{2}}{\sqrt{\lambda t}} \right)
                            + C_2 Y_1\left( \frac{2\sqrt{2}}{\sqrt{\lambda t}} \right) \right), 
 \end{equation}
 where $J_1$ and $Y_1$ are, respectively, the Bessel J- and Y-functions of order 1 (Bessel first
 and second kind functions respectively) and $C_1,C_2$ are general constants. 
 Differentiating for $\lambda \neq 0$ we come to
 \begin{equation}\label{PHIWYN}
  \phi(t) = C_1 \left( \frac{1}{\sqrt{t}}
     J_1\left(\frac{2\sqrt{2}}{\sqrt{\lambda t}}\right) - \frac{\sqrt{2}}
    {t\sqrt{\lambda}} J_0\left(\frac{2\sqrt{2}}{\sqrt{\lambda t}}\right)\right) +
    C_2 \left( \frac{1}{\sqrt{t}}
     Y_1\left(\frac{2\sqrt{2}}{\sqrt{\lambda t}}\right) - \frac{\sqrt{2}}
    {t\sqrt{\lambda}} Y_0\left(\frac{2\sqrt{2}}{\sqrt{\lambda t}}\right)\right).
 \end{equation}
 Recall now that, in small $h>0$ asymptotics, $J_1(h) \sim h/2,\; J_0(h) \sim 1,\;
 Y_1(h) \sim -\frac{2}{\pi h}$ and $Y_0(h) \sim \frac{2}{\pi} \log h,$ see
 e.g. Section 9.4 in Temme (1996). By (\ref{PHIWYN}), for large $t>0$ this readily
 yields $\phi(t) = C_1 o(1/t) + C_2 \Theta(1).$ Likewise, by (\ref{PSIWYN}),
 $\Psi(t) = C_1 \Theta(1) + C_2 \Theta(t).$ Consequently, since
 $\phi \in L_2([1,\infty)),$ we must have $C_2=0.$ In view of (\ref{PSIWYN}) and
 (\ref{WARPOCZ}) this is only possible when
 \begin{equation}\label{WARWWL}
  J_1\left( \frac{2\sqrt{2}}{\sqrt{\lambda}} \right) = 0.
 \end{equation}
 Thus, all eigenvalues of $K$ are positive real numbers satisfying 
 (\ref{WARWWL}). Moreover, they are all simple since, under (\ref{WARWWL}) and
 with $C_2 = 0$ the solution of (\ref{RNIER}) and (\ref{WARPOCZ}) is unique
 up to a multiplicative constant. It remains to check that each $\lambda > 0$ 
 satisfying (\ref{WARWWL}) is an eigenvalue of $K.$ To this end it is enough
 to recall the eigenequation (\ref{RNIEWL}) and observe that its LHS is
 $\lambda \phi(t) = o(1/t)$ and converges to $0$ in large $t$ asymptotics
 and so does the RHS which is asymptotic to $\Psi(t) / t^2 = O(1/t^2)$
 as $t \to \infty.$ This completes the proof of Lemma \ref{SPECK}. $\Box$   

\paragraph{Acknowledgements} 
 The author acknowledges the support from the Polish Minister of Science and Higher 
 Education grant N N201 385234 (2008-2010).
 Special thanks are due to Filip Piekniewski for
 many helpful discussions and for attracting my attention to the experimental
 paper of Egu\'iluz et al. (2005), as well as to Youra Tomilov for
 his helpful comments and remarks.


\paragraph{References}

 \begin{description}
  \item {\sc Albert, R., Barab\'asi, A.-L.} (2002) Statistical mechanics of complex networks, 
          {\it Reviews of modern physics}, {\bf 74}, 47-97.
  \item {\sc Bullmore, E., Sporns, O.} (2009) Complex brain networks: graph theoretical
           analysis of structural and functional systems, {\it Nature Reviews Neuroscience} {\bf 10}, 186-198.
   \item {\sc Cecchi, G.A., Rao, A.R., Centeno, M.V., Baliki, M., Apkarian, A.V., Chialvo, D.R.}
            (2007) Identifying directed links in large scale functional networks: applications to brain
            fMRI, {\it BMC Cell Biology} {\bf 8}, suppl. 1, p.S5.
  \item {\sc Chung, F., Lu, L.} (2006) {\it Complex graphs and networks}, AMS Regional
        Conference Series in Mathematics {\bf 107}, AMS.
  \item {\sc Durrett, R.} (2007) Random Graph Dynamics, {\it Cambridge Series in Statistical and Probabilistic
            Mathematics}, Cambridge University Press. 
  \item {\sc Egu\'iluz, V.M., Chialvo, D.R., Cecchi, G.A., Baliki, M., Apkarian, A.V.}
        (2005) Scale-free brain functional networks, {\it Phys. Rev. Lett.} {\bf 94}, 
        018102.
  \item {\sc Guionnet, A.} (2009) Large Random Matrices: Lectures on Macroscopic Asymptotics,
           \'Ecole d'\'Et\'e de Probabilit\'es de Saint Flour XXXVI - 2006, Lecture Notes in Mathematics.
  \item {\sc Heuvel, van den M.P., Stam, C.J., Boersma, M., Hulshoff Pol, H.E.} (2008) Small
           world and scale-free organization of voxel-based resting-statet functional connectivity in the
           human brain, {\it Neuroimage} {\bf 43}, 528-539. 
  \item {\sc Kato, T.} (1976) {\it Perturbation Theory for Linear Operators}, Grundlehren
        der mathematischen Wissenschaften {\bf 132}, Second Ed., Springer.
  \item {\sc Piersa, J., Piekniewski, F., Schreiber, T.} (2010) Theoretical model for mesoscopic-level scale-free
            self-organization of functional brain networks, submitted.
  \item {\sc Piekniewski, F.} (2007) Emergence of scale-free graphs in dynamical spiking
         neural networks, in: {\it Proc. IEEE International Joint Conference on Neural Networks, Orlando,
         Florida, USA, 2007}, 755-759, IEEE Press.
  \item {\sc Piekniewski, F., Schreiber, T.} (2008) Spontaneous scale-free structure of
         spike flow graphs in recurrent neural networks, {\it Neural Networks} {\bf 21}, 1530-1536.
  \item {\sc Sherrington, D., Kirkpatrick, S.} (1972) Solvable model of a spin glass,
        {\it Phys. Rev. Lett.} {\bf 35}, 1792-1796.
  \item {\sc Shohat, J.A., Tamarkin, J.D.} (1943) {\it The problem of moments}, Math. Survey
        {\bf 1}, AMS, New York.
  \item {\sc Simon, B.} (2005) {\it Trace Ideals and Their Applications}, Math. Surveys and Monographs
        {\bf 120}, AMS, New York. 
  \item {\sc Resnick, S.I.} (1987) {\it Extreme Values, Regular Variation and Point Processes}, 
        Springer-Verlag.
  \item {\sc Salvador, R., Suckling, J., Coleman, M.R., Pickard, J.D., Menon, D., Bullmore, E.} (2005)
          Neurophysiological architecture of functional magnetic resonance images of human brain,
          {\it Cerebral Cortex}, {\bf 15}, 1332-1342.
  \item{\sc Shorack, G.R. and Wellner, J.A.} (1986), {\it Empirical
        Processes with Applications to Statistics}, Wiley, New York. 
  \item {\sc Temme, N.M.} (1996) {\it Special Functions: An Introduction to the Classical
        Functions of Mathematical Physics}, John Wiley \& Sons. 
 \end{description}

\end{document}